# Filtered Ring Derived from Discrete Valuation Ring and Its Properties

M. H. Anjom Shoa*, M. H. Hosseini

University of Birjand, Birjand, Iran
Email: *ajomshoamh@birjand.ac.ir, mhhosseini@birjand.ac.ir

## Abstract

In this paper we show that if $R$ is a discrete valuation ring, then $R$ is a filtered ring. We prove some properties and relation when $R$ is a discrete valuation ring.

## 1. Introduction

In commutative algebra, valuation ring and filtered ring are two most important structures (see [1]-[3]). If $R$ is a discrete valuation ring, then $R$ has many properties that have many usages for example decidability of the theory of modules over commutative valuation domains (see [1]-[3]), Rees valuations, and asymptotic primes of rational powers in Noetherian rings, and lattices (see [4]). We know that filtered ring is also a most important structure since filtered ring is a base for graded ring especially associated graded ring, completion, and some results like on the Andreadakis Johnson filtration of the automorphism group of a free group (see [5]) on the depth of the associated graded ring of a filtration (see [6]). So, as important structures, the relation between these structures is useful for finding some new structure. In this article, we show that we can make a filtration with a valuation. Then we explain some new properties for it. On the other hand, we show this is a strongly filtered ring, then we explain some new properties for it.

## 2. Preliminaries

In this paper the ring $R$ means a commutative ring with unit.

---

*Corresponding author.

**Definition 2.1** A subring $R$ of a field $K$ is called a valuation ring of $K$, if for every $\alpha \in K$, $\alpha \neq 0$, either $\alpha \in R$ or $\alpha^{-1} \in R$.

**Definition 2.2** Let $\Delta$ be a totally ordered abelian group. A valuation $v$ on $R$ with values in $\Delta$ is a mapping $v : R^* \to \Delta$ satisfying:

i) $v(ab) = v(a) + v(b)$;

ii) $v(a+b) \geq \min\{v(a), v(b)\}$.

**Definition 2.3** Let $K$ be field. A discrete valuation on $K$ is a valuation $v : K^* \to Z$ which is surjective.

**Definition 2.4** A fractionary ideal of $R$ is an $R$-submodule $M$ of $K$ such that $aM \subseteq R$, for some $a \in R$, $a \neq 0$.

**Definition 2.5** A fractionary ideal $M$ is called invertible, if there exists another fractionary ideal $N$ such that $MN = R$.

**Proposition 2.1** Let $R$ be a local domain. Every non zero fractionary ideal of $R$ is invertible if and only if $R$ is DVR (see [3]).

**Theorem 2.1** Let $R$ be a Noetherian local domain with unique maximal ideal $m \neq 0$ and $K$ the quotient field of $R$. The following conditions are equivalent.

i) $R$ is a discrete valuation ring;

ii) $R$ is a principal ideal domain;

iii) $m$ is principal;

iv) $R$ is internally closed and every non-zero prime ideal of $R$ is maximal;

v) Every non-zero ideal of $R$ is power of $m$ (see [3]).

**Definition 2.6** Let $R$ be a ring together with a family $\{R_n\}_{n \geq 0}$ of subgroups of $R$ if satisfying the following conditions:

i) $R_0 = R$;

ii) $R_{n+1} \subseteq R_n$ for all $n \geq 0$;

iii) $R_n R_m \subseteq R_{n+m}$ for all $n, m \geq 0$;

Then we say $R$ has a filtration.

**Definition 2.7** Let $R$ be a ring together with a family $\{R_n\}_{n \geq 0}$ of subgroups of $R$ if satisfying the following conditions:

i) $R_0 = R$;

ii) $R_{n+1} \subseteq R_n$ for all $n \geq 0$;

iii) $R_n R_m = R_{n+m}$ for all $n, m \geq 0$;

Then we say $R$ has a strong filtration.

**Example 2.1** Let $I$ be an ideal of $R$, then $R_i = I^i$ is a filtration that is called $I$ adic filtration ring.

**Definition 2.8** Let $R$ be a filtered ring. A filtered $R$-module $M$ is an $R$-module together with family $\{M_n\}_{n \geq 0}$ of subgroup $M$ of satisfying:

1. $M_0 = M$;

2. $M_{n+1} \subseteq M_n$ for all $n \geq 0$;

3. $R_n M_m \subseteq M_{n+m}$ for all $n, m \geq 0$.

Then we say $M$ has a filtration.

**Definition 2.9** A map $f : M \to N$ is called a homomorphism of filtered modules, if: i) $f$ is $R$-module an homomorphism and ii) $f(M_n) \subseteq N_n$ for all $n \geq 0$.

**Definition 2.10** A graded ring $R$ is a ring, which can expressed as a direct sum of subgroup $\{R_n\}_{n \geq 0}$ i.e. $R = \oplus_{n \geq 0} R_n$ such that $R_n R_m \subseteq R_{n+m}$ for all $n, m \geq 0$

**Definition 2.11** Let $R$ be a graded ring. An $R$-module $M$ is called a graded $R$-module, if $M$ can be expressed as a direct sum of subgroups $\{M_n\}_{n \geq 0}$ i.e. $M = \oplus_{n \geq 0} M_n$ such that $R_n M_m \subseteq M_{n+m}$ for all $n, m \geq 0$.

**Definition 2.12** Let $M$ and $N$ be graded modules over a graded ring $R$. A map $f : M \to N$ is called homomorphism of graded modules if: i) $f$ is $R$-module an homomorphism and ii) $f(M_n) \subseteq N_n$ for all $n \geq 0$.

**Definition 2.13** Let $R$ be a filtered ring with filtration $\{R_n\}_{n \geq 0}$. Let $gr_n(R) = R_n/R_{n+1}$, and $gr(R) = \oplus_{n \geq 0} gr_n(R)$. Then $gr(R)$ has a natural multiplication induced from $R$ given

$$(a + R_{n+1})(b + R_{m+1}) = ab + R_{m+n+1}$$

where $a \in R_n, b \in R_m$. This makes $R$ in to a graded ring. This ring is called the associated graded ring of $R$.

**Definition 2.14** Let $M$ be a filtered $R$-module over a filtered ring $R$ with filtration $\{M_n\}_{n \geq 0}$ and $\{R_n\}_{n \geq 0}$ respectively. Let $gr_n(M) = M_n/M_{n+1}$, and $gr(M) = \oplus_{n \geq 0} gr_n(M)$. Then $gr(M)$ has a natural $gr(R)$-module structure given by $a + R_{n+1}(x + M_{m+1}) = ax + M_{m+n+1}$, where $a \in R_n, x \in M_m$.

## 3. Filtered Ring Derived from Discrete Valuation Ring and Its Properties

In this section we proved that, if $R$ is a discrete valuation ring, then $R$ is a filtered ring. And we prove some properties for $R$.

Let $K$ be a field which $R$ be a domain and a discrete valuation ring (DVR) for $K$. The map $\upsilon : K^* \to Z^+$ is valuation of $R$.

**Lemma 3.1** By above definition, the set $R_n = \{\alpha \in K | \nu(\alpha) \geq n, n \in \mathbb{Z}\}$ is an ideal of $R$.

**Proof.** (see [3])

**Theorem 3.1** If $R$ is a discrete valuation ring with valuation $\nu : K^* \to Z^+$. Then $R$ is a filtered ring with filtration defined by

$$R_n = \{\alpha \in K | \nu(\alpha) \geq n, n \in \mathbb{Z}\}$$

where $R_0 = R$.

**Proof.** By definition of valuation ring, it is obvious that $R_0 = R$. For the second condition for filtration ring we have $\forall \alpha \in R_{n+1} \Rightarrow \nu(\alpha) \geq n+1 > n \Rightarrow \nu(\alpha) > n \Rightarrow \alpha \in R_n$, So we have $R_{n+1} \subseteq R_n$.

For the third condition, we have for every $R_n$ and $R_m$ without losing generality. Since $R_n$ and $R_m$ are ideals of $R$ so

$$R_n R_m = \{\sum a_i b_i | a_i \in R_n - b_i \in R_m\}$$

is an ideal of $R$.

Now let $c \in R_n R_m$ then $c = \sum_{i \in I} a_i b_i$ for $a_i \in R_n$ and $b_i \in R_m$. Thus

$$\nu(c) = \nu\left(\sum_{i \in I} a_i b_i\right) = \min\{\nu(a_i b_i)\}_{i \in I} = \min\{\nu(a_i) + \nu(b_i)\}_{i \in I} \geq \left(\min\{\nu(a_i)\}_{i \in I}\right) + \left(\min\{\nu(b_i)\}_{i \in I}\right) \geq n + m,$$

Consequently we have $\nu(c) \geq n + m \Rightarrow c \in R_{n+m}$ hence $R_n R_m \subseteq R_{n+m}$. Therefore $R$ is a filtered ring.

**Proposition 3.1** Let $R$ be a local domain. If every non-zero fractionary ideal of $R$ invertible, then $R$ is filtered ring.

**Proof.** By proposition 2.1 $R$ is DVR then by theorem 3.1 $R$ is filtered ring.

**Proposition 3.2** Let $R$ be a Noetherian local domain with unique maximal ideal $m \neq 0$ and $K$ the quotient field of $R$. Then $R$ is filtered ring if one of following conditions is held
  i) $R$ is a principal ideal domain;
  ii) $m$ is principal;
  iii) $R$ is integrally closed and every non-zero prime ideal of $R$ is maximal.

**Proof.** It follows from theorem (3.1) and theorem (2.1).

**Definition 3.1** Let $R$ be a ring, and let $(S,+)$ be a totally ordered cancellative semigroup having identity $0$. A function $f : R \to S \cup \{\infty\}$ is a filtration if $f(1) = 0$, $f(0) = \infty$ and for all $x, y \in R$,

i) $f(x+y) \geq \min\{f(x), f(y)\}$, and

ii) $f(xy) \geq f(x) + f(y)$, then $f$ is called a filtration.

For this filtration we have

1) $\{A_g : g \in S^+\}$ the set of ideals;

2) $A_+ = \{x \in R : f(x) > 0\}$;

3) $(f)^g = \{x \in R : \exists n > 0 \text{ such that } g \leq f(x^n)\}$;

4) $(f)_g = \{x \in R : \exists n > 0 \text{ such that } g = f(x^n)\}$.

**Lemma 3.2** Let $f : R \to S \cup \{\infty\}$ be a filtration and let $g, h \in S$. Then:

i) $(f)_0 = \sqrt{A_+}$;

ii) $(f)^\infty = 0$;

iii) $(f)_g \subseteq (f)^g$;

iv) if $g \leq h$, then $(f)_h \subseteq (f)_g$ and $(f)^h \subseteq (f)^g$.

**Proof.** See lemma 3.3 of [7].

**Proposition 3.3** If $R$ be a discrete valuation ring, then there exists a totally ordered cancellative semigroup $S$, and $f : R \to S \cup \{\infty\}$ such that:

i) $(f)_0 = \sqrt{A_+} = \sqrt{m}$;

ii) $(f)^\infty = 0$;

iii) $(f)_g \subseteq (f)^g$;

iv) if $g \leq h$, then $(f)_h \subseteq (f)_g$, and $(f)^h \subseteq (f)^g$.

**Proof.** By theorem 3.1 there exists a filtration for $R$, then by lemma 2.1 we have the all above conditions.

**Proposition 3.4** Let $R$ be a filtered ring, $M$, $N$ filtered $R$-modules, and $f : M \to N$ homomorphism of filtered $R$-modules. If the induced map $gr(f) : gr(M) \to gr(N)$ is injective, then $f$ is injective provided $\bigcap_{n=0}^{\infty} M_n = (0)$. (see [3])

**Corollary 3.1** Let $R$ be a valuation ring, $M$, $N$ filtered $R$-modules, and $f : M \to N$ homomorphism of filtered $R$-modules. If the induced map $gr(f) : gr(M) \to gr(N)$ is injective, then $f$ is injective provided $\bigcap_{n=0}^{\infty} M_n = (0)$.

**Proposition 3.5** If $R$ is a discrete valuation ring with valuation $v : K^* \to Z^+$, Then $R$ is a strongly filtered ring with filtration defined by

$$R_n = \{\alpha \in K | v(\alpha) \geq n, n \in \mathsf{Z}\}$$

where $R_0 = R$.

**Proof.** By theorem 3.1 $R$ is a filtered ring. Now we show $R_n R_m = R_{n+m}$ for all $\alpha \in R_{n+m}$. Since $n + m \geq n, m$ so

$$R_{m+n} \subseteq R_n \text{ and } R_{n+m} \subseteq R_m.$$

Consequently $\alpha \in R_n$, and $\alpha \in R_m$. Therefore $R_{n+m} \subseteq R_n R_m$.

**Proposition 3.6** Let $R$ be a discrete valuation ring, and $f : R \to S \cup \{\infty\}$. If $x \in R$ and $f(x) > 0$, then $(f)^{f(x)}$ is smallest prime ideal in $Spec_f(R)$ which contains $x$, and $(f)_{f(x)}$ is largest prime ideal in $Spec_f(R)$ which does not contains $x$.

**Proof.** By proposition 3.5 $R$ is strongly filtered ring, then by proposition 4.2. of [7]-[9] we have If $x \in R$ and $f(x) > 0$, then $(f)^{f(x)}$ is smallest prime ideal in $Spec_f(R)$ such that contains $x$, and $(f)_{f(x)}$ is the largest prime ideal in $Spec_f(R)$ such that does not contains $x$.

**Remark 3.1** Given a strong filtration $f$ on a ring $R$, we say that a prime $P$ in $Spec_f(R)$ is branched in $Spec_f(R)$, if $P$ cannot be written as union of prime ideals in $Spec_f(R)$ such that properly contained in $P$.

**Corollary 3.2** Let $R$ be a discrete valuation ring and $f : R \to S \cup \{\infty\}$. Then a prime ideal $P$ in $Spec_f(R)$ is branched in $Spec_f(R)$ if and only if $P = (f)^g$ for some $g \in S^+$.

**Proof.** By proposition 3.5 $R$ is strongly filtered ring, then by proposition 4.5. of [7] a prime ideal $P$ in $Spec_f(R)$ is branched in $Spec_f(R)$, if and only if, $P = (f)^g$ for some $g \in S^+$.